%% file: Disque-CRAS-15-03-31.tex
\documentclass[10pt]{article}
\usepackage[T1]{fontenc}
\usepackage{amssymb,amsmath,amsfonts}
\usepackage{dsfont}
\usepackage{amsthm}
\usepackage{graphicx,color,subfig}
\usepackage{xspace}
\usepackage{txfonts}
\usepackage{mathrsfs}
\usepackage{amscd}
\usepackage{amssymb}

\usepackage{mathrsfs}
\usepackage[colorlinks=true]{hyper ref}

\newtheorem{theorem}{Theorem}[section]
\newtheorem{corollary}[theorem]{Corollary}
\newtheorem{definition}[theorem]{Definition}

\newtheorem{lemma}[theorem]{Lemma}

\newtheorem{proposition}[theorem]{Proposition}
\newtheorem{remark}[theorem]{Remark}

\newcommand{\ovl}[1]{\overline{#1}}
\newcommand{\W}{\ensuremath{\mathbb W}}
\renewcommand{\d}{\ensuremath{\partial}}
\newcommand{\scrU}{\ensuremath{\mathscr U}}

\DeclareMathOperator{\loc}{loc}

\input{def.tex}

\numberwithin{equation}{section}

\makeatletter
\def\keywords{
    \vspace{1ex}
    \noindent
    \if@twocolumn
      \small{\bf  Keywords}\/---$\!$    \else
      \begin{center}\small\ {\bf Keywords}\end{center}\quotation\small
    \fi}
\def\endkeywords{\vspace{0.6em}\par\if@twocolumn\else\endquotation\fi
    \normalsize\rm}
\makeatother

\begin{document}

\title{Delocalization of quasimodes on the disk\footnote{
\textbf{Acknowledgement.} NA and ML are partially supported by the Agence Nationale de la Recherche under grant GERASIC ANR-13-BS01-0007-01.
FM is partially supported by grants MTM2013-41780-P (MEC) and ERC Starting Grant 277778.
}
}

\author{
Nalini Anantharaman\footnote{Universit\'{e} de Strasbourg, IRMA, 7 rue Ren\'e-Descartes, 67084 Strasbourg Cedex, France. e-mail: {\tt anantharaman@unistra.fr}},
Matthieu L\'eautaud\footnote{Universit\'e Paris Diderot, Institut de Math\'ematiques de Jussieu-Paris Rive Gauche, UMR 7586, B\^atiment Sophie Germain, 75205 Paris Cedex 13 France. e-mail: {\tt leautaud@math.univ-paris-diderot.fr}},
Fabricio Maci\`a\footnote{Universidad Polit\'{e}cnica de Madrid. DCAIN, ETSI Navales. Avda. Arco de la
Victoria s/n. 28040 Madrid, Spain. e-mail: {\tt Fabricio.Macia@upm.es}}
}

\maketitle

\begin{abstract}
This note deals with semiclassical measures associated to {(sufficiently accurate)} quasimodes $(u_h)$ for the Laplace-Dirichlet operator on the disk. In this time-independent set-up, we simplify the statements of~\cite{ALM:disk} and their proofs. We describe the restriction of semiclassical measures to every invariant torus in terms of two-microlocal measures. As corollaries, we show regularity and delocalization properties for limit measures of $|u_h|^2 dx$: these are absolutely continuous in the interior of the disk and charge every open set intersecting the boundary.
\end{abstract}

\begin{abstract}
\textbf{D\'elocalisation des quasimodes sur le disque.}
Dans cette note, on s'int\'eresse aux mesures semiclassiques associ\'ees aux quasimodes (d'ordre suffisamment \'elev\'e) $(u_h)$ du laplacien de Dirichlet sur le disque. Dans ce contexte stationnaire, les r\'esultats obtenus dans~\cite{ALM:disk} et leurs preuves sont simplifi\'es. On d\'ecrit la restriction de ces mesures \`a chaque tore invariant au moyen de mesures deux-microlocales. En corollaire, on montre des propri\'et\'es de r\'egularit\'e et de d\'elocalisation des mesures limites des $|u_h|^2 dx$~: celles-ci sont absolument continues \`a l'int\'erieur du disque et chargent tout ouvert qui touche le bord.
\end{abstract}

%
%

\section{Introduction}
We consider the unit disk
$
\ID = \{z=(x,y) \in \IR^2 , |z|^2 = x^2 + y^2 <1\} \subset \IR^2 ,
$
and study quasimodes for the
euclidean Laplacian $\Delta$ endowed with Dirichlet boundary conditions:
\begin{eqnarray}
\label{e:S} 
 \left( - h^2\Delta + h^2 V - E_0^2 \right)u_h = r_h,  \quad \text{ in }\ID,
 \quad  u_h|_{\d \ID} = 0  , \quad \|u_h\|_{L^2(\ID)}  = 1 .
\end{eqnarray}
where $V = V(z)$ is a bounded potential and $E_0 >0$ a fixed energy level (say $E_0 =1$). Here, $h>0$, $h \to 0$ is a semiclassical parameter and the remainder $r_h$ satisfies some boundedness/smallness assumptions in $L^2(\ID)$: 
\begin{definition}
Let $(s_h)$ be a family of positive real numbers indexed by $h\in (0,1)$. We say that $(u_h)_{h>0}$ is a family of $O(s_h)$ (resp. $o(s_h)$) quasimodes if $u_h$ satisfies~\eqref{e:S} with $\|r_h\|_{L^2(\ID)} =O(s_h)$ (resp. $\|r_h\|_{L^2(\ID)} =o(s_h)$) as $h \to 0^+$.
\end{definition}
The aim of this note is to prove delocalization properties for sufficiently accurate quasimodes, namely $O(h^2)$ or $o(h^2)$ quasimodes.
This type of result can be deduced from similar properties for solutions to the semiclassical Schr\"odinger evolution equation 
$
\frac{h}{i} \d_t w_h = h^2\left(-\Delta+V\right) w_h
$.
If $(u_h)$ solves~\eqref{e:S}, then the solution $w_h(t)$ of the evolution equation with $w_h|_{t=0} =u_h$ satisfies $\| w_h(t) - e^{it E_0/h}u_h \|_{L^2(\ID)}= O(\frac{t}{h} \|r_h\|_{L^2(\ID)})$. Hence, properties of $w_h(t)$ over a time interval $[0,\tau_h]$ can be translated into properties of $u_h$ if $\|r_h\|_{L^2(\ID)} \sim \frac{h}{\tau_h}$. The article~\cite{ALM:disk} deals with properties of solutions of the time dependent Schr\"odinger equation on time intervals of length $1/h$; the results of~\cite{ALM:disk} can therefore be transferred into properties for quasimodes of order $h^2$ (see also~\cite[Remark~2.5]{ALM:disk}).
Note that, although all the results we present are special cases of those in~\cite{ALM:disk}, considering stationary solutions allows to simplify the statements and the proofs significantly. This is the motivation of the present note.

Examples of quasimodes are provided by high-energy eigenfunctions of $\left( - \Delta+V \right)$ or clusters of eigenfunctions: denote by $(\psi_j)$ a Hilbert basis of $L^2(\ID)$ consisting of orthonormal eigenfunctions of the operator $-\Delta + V$ with Dirichlet conditions, associated to the eigenvalues $\lambda_j \to + \infty$. Then the family
$$
v^{\lambda} = \sum_{\lambda_j \in [\lambda-R(\lambda), \lambda +R(\lambda)]} v_j \psi_j ,  \quad \text{ with }  \sum_{\lambda_j \in [\lambda-R(\lambda), \lambda+R(\lambda)]} |v_j|^2 = 1, \quad \lambda \to + \infty ,
$$ 
is a family of $O(s_h)$ quasimodes for $s_h = R(\lambda)/\lambda$ and $h = E_0 \lambda^{-1/2} \to 0$.

A major issue in mathematical quantum mechanics is to describe the possible localization -- or delocalization -- properties of solutions of the stationary Schr\"odinger equation~\eqref{e:S}. 
Here, the main object of our study is the probability density $|u_h(z)|^2 dz$; given $\Omega \subset \ID$, the quantity $\int_{\Omega}|u_h(z)|^2 dz$ represents the probability of finding a quantum particle in the set $\Omega$. 
More precisely, given a sequence $h = h_n \to 0^+$, we aim at describing the asymptotic properties of the probability densities $|u_h(z)|^2 dz = |u_{h_n}(z)|^2 dz$. 
After possibly extracting a subsequence, we have the convergence $|u_h(z)|^2 dz \rightharpoonup \nu(dz)$ in $\mathcal{D}'(\R^2)$, where $\nu$ is a nonnegative Radon measure describing the asymptotic mass repartition of the sequence of quasimodes $(u_h)$. One of the goals of this paper is to understand how the fact that $(u_h)$ solves~\eqref{e:S} influences the structure of the associated measure $\nu$. 

Another interesting quantity is the mass left by a quasimode at the boundary: a well-known {\em hidden regularity} result (see e.g.~\cite[Lemma~2.1]{GerLeich93}) states that $h\d_n u_h|_{\d \ID} $ forms a bounded family of $L^2(\d \ID)$ for any family of $O(1)$ quasimodes $(u_h)$. Hence one may also be interested in studying the asymptotic repartition of the densities $|h\d_n u_h|_{\d \ID}|^2 dS(z)$, where $dS$ denotes the Lebesgue measure on the circle $\d \ID$. After extracting a subsequence, one has $|h\d_n u_h|_{\d \ID}|^2 dS(z)\rightharpoonup \nu^\d(dz)$ where $\nu^\d$ is a measure on the boundary $\d \ID$.

Theorem~\ref{t:precise} in Section \ref{section: main} describes precisely the properties of semiclassical measures: these are lifts of the limit measures $\nu$ (described above) to the phase space of classical dynamics. The theorem deals with $O(h^2)$ or $o(h^2)$ quasimodes and it yields in particular the following three corollaries.
\begin{corollary}
\label{t:example}Let $(u_h)$ be a sequence of $O(h^2)$ quasimodes. 

(i) For every weak-$\ast$ limit $\nu(dz)$ of the sequence $|u_h(z)|^2 dz$, the restriction $\nu|_{\ID}$ is absolutely continuous.

(ii) Any weak-$\ast$ limit $\nu^\d(dz)$ of the sequence $|h\d_n u_h|_{\d \ID}|^2 dS(z)$ is absolutely continuous (with respect to $dS$). 
\end{corollary}
This result shows that the weak-$\ast$ accumulation points of the
densities $|u_h(z)|^2 dz$ possess some regularity in the interior
of the disk (note that
it is easy to exhibit sequences of quasimodes that
concentrate singularly on the boundary, the so-called
whispering-gallery modes, having for limit measure $\nu(dz) = (2\pi)^{-1} \delta_{\d \ID}$). 
Remark that a family of $O(h^2)$ quasimodes for $- h^2\Delta + h^2 V - E_0^2$ is a family of $O(h^2)$ quasimodes for $- h^2\Delta - E_0^2$. As Corollary~\ref{t:example} applies for $O(h^2)$ quasimodes, no regularity is needed for the potential $V$ and the result also holds under the assumption $V \in \mathcal{L}(L^2(\ID))$.

Such a regularity result is also known to hold on flat tori \cite{MaciaTorus,AnantharamanMaciaTore} and more generally in the case of strictly convex/concave completely integrable systems (without boundary)~\cite{AnantharamanKMacia}. 
On the sphere $\mathbb{S}^d$, on which the geodesic flow is still completely integrable, the situation is radically different, for it is known that every measure that is invariant under the geodesic flow (in particular, the uniform measure on an equator) is a semiclassical measure.

Note that it is proved in~\cite[Theorem~1.3]{AnantharamanKMacia} that the scale $h^2$ is the critical delocalization scale for quasimodes on non degenerate completely integrable systems: $O(s_h)$ quasimodes with $s_h \gg h^2$ can have as a semiclassical measure every invariant measure of the geodesic flow. {In that reference it is also shown that the size $h^2$ of the potential is also critical: it is possible to give an example of a potential $V$ such that for any $\eps>0$ there exists a sequence of $O(h^\infty)$ quasimodes $(u_h)$ for the operator $- h^2\Delta + f(h)V$,  with $f(h)=O(h^{2-\eps})$, such that $|u_h(z)|^2 dz$ concentrates singularly on a classical trajectory.}

Another corollary of Theorem~\ref{t:precise} is the following result:
\begin{corollary}
\label{t:obs-i} 
(i) Let $\Omega\subset \ovl{\ID}$ be an open set such
that $\Omega\cap \partial \ID \not= \emptyset$, and $V \in C^\infty(\overline \ID ; \R)$. Then, there exist $C(\Omega)>0$ such that 
for any sequence $(u_h)$ of $o(h^2)$ quasimodes, for any weak-$\ast$ limit $\nu(dz)$ of the sequence $|u_h(z)|^2 dz$, we have $\nu (\Omega) \geq C(\Omega)$. 

(ii) Let $\Gamma \subset \d \ID$ be any nonempty open set, and $V \in C^\infty(\overline \ID; \R)$. Then, there exist $C(\Gamma)>0$ such that for any sequence $(u_h)$ of $o(h^2)$ quasimodes, for any weak-$\ast$ limit $\nu^\d(dz)$ of the sequence $|h \d_n u_h(z)|_{\d \ID}|^2 dS(z)$, we have $\nu^\d(\Gamma) \geq C(\Gamma)$. 
\end{corollary}

Points (i) and (ii) of Corollary \ref{t:obs-i} are equivalent (after {\em reductio ad absurdum} and the use of unique continuation for eigenfunctions of the operator $-\Delta + V$) to the following resolvent estimates:

\begin{corollary}
\label{t:resolv-estim}
(i) Let $\Omega\subset \ovl{\ID}$ be an open set such
that $\Omega\cap \partial \ID \not= \emptyset$, and $V \in C^\infty(\overline \ID; \R)$. Then, there exist $C_0, C_1>0$ such that for any $\lambda \in \R$, for any $u \in H^2 \cap H^1_0(\ID)$ we have
$$
\|u\|_{L^2(\ID)} \leq C_0\|(-\Delta + V - \lambda)u\|_{L^2(\ID)} + C_1\|u\|_{L^2(\Omega)}
$$

(ii) Let $\Gamma \subset \d \ID$ be any nonempty open set, and $V \in C^\infty(\overline \ID; \R)$.
Then, there exist $C_0, C_1>0$ such that for any $\lambda \in \R$, for any $u \in H^3 \cap H^1_0(\ID)$ such that $\Delta u|_{\d \ID} = 0$, we have
$$
\|u\|_{H^1_0(\ID)} \leq C_0\|(-\Delta + V - \lambda)u\|_{H^1_0(\ID)} + C_1\|\d_n u|_{\d \ID}\|_{L^2(\Gamma)}
$$
\end{corollary}
Roughly speaking, this means that any set $\Omega$ touching $\d \ID$ (resp. any subset $\Gamma$ of $\d \ID$) observes all quantum particles trapped in the disk. Because of the whispering gallery phenomenon, the condition that $\Omega$ touch the boundary is necessary for property (i) to hold. 
This reflects the fact that any solution has to leave positive mass on any set $\Omega$ touching the boundary $\d \ID$ (resp. any subset $\Gamma$ of $\d \ID$).
In the present very particular geometry, this improves the general bound~\cite{LR:95} (given by the tunnelling effect) where $C_0, C_1$ have to be replaced by $C e^{C\lambda}$ for some $C>0$.
Resolvent estimates such as those of Corollary~\ref{t:resolv-estim} are known to imply observability/controllability results for the evolution Schr\"odinger equation in sufficiently large time~\cite{BurqZworskiBlackBox04,RTTT:05, Miller:05}. 

It is known that the resolvent estimates of Corollary~\ref{t:resolv-estim} hold in a general domain (in an improved form, with $C_0$ replaced by $C_0(1+|\lambda|)^{-1}$) under the stronger assumption that all trajectories of the billiard enter the observation region $\Omega$ or $\Gamma$ in finite time~\cite{Leb:92,BurqZworskiBlackBox04,RTTT:05, Miller:05}. 
There are other situations in which this strong geometric control condition is not necessary. This is the case for the torus, for (i) is satisfied as soon as $\Omega \neq \emptyset$~\cite{JaffardPlaques,MaciaDispersion,BZ:12,Kom:92,AnantharamanMaciaTore}. The boundary resolvent estimate of Corollary~\ref{t:resolv-estim} also holds in the square if and only if the observation
region $\Gamma$ contains both a horizontal and a vertical nonempty
segments~\cite{RTTT:05}. On the other hand, on the sphere, it is {\em necessary} that $\overline{\Omega}$ meets all geodesics for an observation inequality like that of Corollary~\ref{t:resolv-estim}
to hold.

\begin{remark}
(i) Arguments developed in \cite{AnantharamanMaciaTore} show that
Corollaries~\ref{t:obs-i} and~\ref{t:resolv-estim} (as well as Theorem~\ref{t:precise} below) also hold for $V\in C^0\left( \ovl{\ID}; \R \right)$ or even in the case where $V$ is continuous outside a set of zero measure. 

(ii) Corollary~\ref{t:resolv-estim} directly yields a polynomial decay rate for
the energy of the (internally) damped wave equation on the disk if the damping touches the boundary~\cite[Theorem~2.3]{AnL}.
\end{remark}

\section{The billiard flow in the disk, and associated Action-Angle coordinates\label{s:billiard}}
Semiclassical analysis provides a connection between
quasimodes and the billiard on the underlying phase space. 
Let us clarify what we mean by ``billiard flow'' in the disk.
 We first define the symmetry with respect to the line tangent to the circle $\partial \ID$ at $z \in \partial \ID$ by $\sigma_z(\xi) = \xi - 2z\cdot \xi$ for $z \in \partial \ID$.
Then, we work on the quotient space
$\W = \ovl{\ID}\times \R^2 / \sim$ where $(z, \xi) \sim
(z,\sigma_z(\xi)) \text{ for } |z| = 1$. 
We denote by $\pi$ the canonical projection $\ovl{\ID}\times \R^2
\to \W$ which maps a point $(z, \xi)$ to its equivalence class
modulo $\sim$. Note that $\pi$ is one-one on $\ID\times \R^2$, so
that $\ID\times \R^2$ may be seen as a subset of $\W$.
A function $a \in C^0(\W)$ can be identified with the function
$\tilde{a} = a \circ \pi \in C^0(\ovl{\ID}\times \R^2)$ satisfying
$\tilde{a} (z, \xi) = \tilde{a}\circ \sigma_z (\xi)$ for $(z,
\xi) \in \d \ID\times \R^2$.
The billiard flow $(\phi^\tau)_{\tau\in \R}$ on $\W$ is the
(uniquely defined) action of $\R$ on $\W$ such that the map
$(\tau, z, \xi)\mapsto \phi^\tau(z, \xi)$ is continuous on
$\IR\times \W$, satisfies
$\phi^{\tau+\tau'}=\phi^\tau\circ\phi^{\tau'}$, and such that
$\phi^\tau(z, \xi)=(z+\tau\xi, \xi)$ 
whenever $z\in\ID$ and $z+\tau\xi\in\ID$.

\medskip
In order to understand how the completely integrable dynamics of
the flow $\phi^\tau$ influences the structure of Wigner measures,
we need to introduce coordinates adapted to these dynamics. We
denote by $\Phi :(s,\theta,E,J)\mapsto(x,y,\xi_{x},\xi_{y})$
the set of ``action-angle'' coordinates for the billiard flow, defined by:
\begin{equation*}
\begin{cases}x=\frac{J}{E}\cos\theta-s\sin\theta,\\
y=\frac{J}{E}\sin\theta+s\cos\theta,\\
\xi_x=-E\sin\theta,\\
\xi_y=E\cos\theta.
 \end{cases}
 \Longleftrightarrow \ \ 
\begin{cases}
E=\sqrt{\xi_{x}^{2}+\xi_{y}^{2}},\mbox{ (velocity)}\\
J=x\xi_{y}-y\xi_{x} = z \cdot \xi^\perp,\mbox{ (angular momentum)}\\
\theta=-\arctan\left(  \frac{\xi_{x}}{\xi_{y}}\right),
\mbox{ (angle of $\xi$ with the vertical)}\\
s=-x\sin\theta+y\cos\theta,\mbox{ (abscissa of }(x,y)\mbox{ along
the line $\left(\frac{J}{E}\cos\theta,
\frac{J}{E}\sin\theta\right) + \R\xi$ ).}
\end{cases}
\]
Above, we have denoted $\xi^\perp = (\xi_{y} , -\xi_{x})$, where $\xi= (\xi_x , \xi_y)$.
Note that the velocity $E$  and the angular
momentum $J$ are preserved along the free transport flow in $\R^2
\times \R^2$, but also along $\phi^\tau$; the variables $s$ and
$\theta$ play the role of ``angle'' coordinates. We call
$\alpha=-\arcsin\left(\frac{J}E\right) = -\arcsin \left(\frac{x\xi_y - y \xi_x}{|\xi|}\right)$ the angle that a billiard
trajectory makes with the normal to the circle, when it hits the
boundary. The quantity $\alpha$
is preserved by the billiard flow.

We set $X_J  =  z^\perp \cdot \d_z + \xi^\perp \cdot \d_\xi$ and $X_E = \frac{\xi}{|\xi|} \d_z$ to be the Hamiltonian vector fields associated to $J(z,\xi)$ and $E(z,\xi)$, respectively. Note that $R^\tau$, the flow of $X_J$, is given by $R^\tau(z, \xi)= (R(\tau)z, R(\tau)\xi)$, where $R(\tau)$ is the rotation matrix of angle $\tau$. Let us denote $T_{(E, J)}$ the level sets of the pair $(E, J)$,
namely
\begin{equation*}
T_{(E, J)} = \{(z,\xi) \in \ovl{\ID} \times \R^2 \;:\; (|\xi|, z
\cdot \xi^\perp) =  (E,J) \}.
\end{equation*}
For $E\not=0$ let us denote $\lambda_{E, J}$ the probability
measure on $T_{(E, J)}$ that is both invariant under the billiard
flow and invariant under rotations. In the coordinates $(s,
\theta, E, J)$, we have
$$
\lambda_{E, J}(ds, d\theta)  = c(E,J)ds  d\theta , \quad c(E,J) =
\left(\int_{T(E,J)}ds d\theta \right)^{-1}>0 .
$$
Note that for $E\not= 0$ and $\alpha\in \pi \IQ$ the billiard flow
is periodic on $T_{(E, J)}$ whereas $\alpha\not\in \pi \IQ$
corresponds to trajectories that hit the boundary on a dense set.
More precisely, if $\alpha\not\in \pi \IQ$ then the billiard flow
restricted to $T_{(E, J)}$ has a unique invariant probability
measure, namely $\lambda_{E, J}$.
For each $\alpha_0\in \pi \IQ\cap (-\pi/2, \pi/2)$ we 
define 
$$
\cI_{\alpha_0}=\{(s, \theta, E, J) \in \Phi^{-1}(\ovl{\ID} \times \R^2) , J=-\sin\alpha_0 E\} = \{\alpha = \alpha_0\},
$$ which is the union of all the lagrangian manifolds $T_{(E, J)}$ with $J=-\sin\alpha_0 E$.
The billiard flow $\phi^\tau$ is periodic on $\cI_{\alpha_0}$; hence, given a function $a : \ID \times \R^2 \to \C$, we may define $\la a\ra_{\alpha_0} : \cI_{\alpha_0} \to \C$ its average along the orbits of
$\phi^\tau$ on the set $\cI_{\alpha_0}$.
In the coordinates $(s, \theta, E, J)$, this function only depends on $\theta$ and $E$.

In the following, we need to perform semiclassical analysis in the variables $(s,\theta,J,E)$ instead of $(z,\xi)$ and hence to quantize the symplectic change of variables $\Phi$. 
\begin{lemma}
\label{p:FIO}
There exist a Fourier Integral Operator $\mathscr U$ satisfying

 (i) $\Op_h(a(z,\xi)) = \mathscr U^{\ast} \Op_h(a \circ \Phi(s,\theta,J,E))  \mathscr U + O(h)$ for any $a \in C^\infty_c(T^*\R^2)$ supported away from $\xi =0$.

 (ii) The operator $\mathscr U$ is unitary from $L^{2}%
(\mathbb{R}^{2})$ to $L^{2}\left(  \R\times\R/2\pi\Z\right) $:
$\mathscr U^{\ast}\mathscr U=I$.

(iii) For $f\in C_{c}^{\infty}(\R^{2})$, we have $\partial^{2}_{s}
\mathscr U f=\mathscr U \Delta f$, i.e. $-h^{2}\mathscr U\Delta\mathscr U^{\ast}=-h^{2}\partial _{s}^{2}$.

\end{lemma}

\section{Semiclassical measures and the structure theorem}
\label{section: main}
Following~\cite{GerLeich93}, we extend the problem from $\ID$ to $\R^2$: starting with $u_h \in H^2\cap H^1_0(\ID)$ we extend this function to $\R^2$ by the value $0$ outside $\ID$ (we still denote by $u_h$ the extended function). The extended function $u_h$ satisfies $u_h \in H^1(\R^2)$, as well as
\begin{eqnarray}
\label{e:SR2} 
 \left( - h^2\Delta + h^2 V - E_0^2 \right)u_h = r_h + h^2 \d_n u_h|_{\d \ID} \otimes \delta_{\d \ID} ,  \quad \text{ in }\R^2  , \quad \|u_h\|_{L^2(\R^2)}  = 1 .
\end{eqnarray}
The \emph{semiclassical Wigner distribution} associated to
$u_{h}$ (at scale $h$) is a distribution on the cotangent bundle
$T^{\ast}\IR^2=\IR^2_z\times \IR^2_\xi$, defined by
\begin{equation*}
W_h : a\mapsto
 \left\langle
u_{h},\Op_h(a(z, \xi))u_{h}\right\rangle
_{L^{2}(\R^2)},\qquad \mbox{ for all }a\in
C_{c}^{\infty}(T^{\ast}\IR^2),\label{e:defW}
\end{equation*}
where $\Op_h$ denotes the standard semiclassical quantization.
After possibly extracting a subsequence, we have, 
\begin{equation}\label{e:wlimite}
W_h(a)\to \mu(a), \quad \mbox{ as }h\to0 \quad \text{for all } a\in
C_c^\infty\left(  T^{\ast}\IR^2\right),
\end{equation} 
where $\mu$ is a nonnegative measure on $T^*\R^2$ called the semiclassical measure associated to the subsequence $(u_h)$.
Our main goal is to describe as precisely as possible the semiclassical measures $\mu$ associated to quasimodes.
It follows from~\cite{GerLeich93} that the limit $\mu$ in \eqref{e:wlimite} has the
following properties (on any convex domain):

$(i)$ If $(u_h)$ is a family of $o(1)$ quasimodes, then $\mu$ is a nonnegative probability measure supported in $S^*_{E_0}\ovl\ID = \{(x, \xi) \in T^* \R^2, x \in\ovl{\ID},  |\xi|=E_0\}$.

$(ii)$ If $(u_h)$ is a family of $o(h)$ quasimodes, then we have $\int_{\ovl{\ID}\times
\R^2\times\R}\xi\cdot\partial_z a\,\, \mu (dz,d\xi)=0$ for every smooth $a$ such
that $a(z, \xi )=a(z, \sigma_z(\xi) )$ for $|z|=1$.
Equivalently,
$\int_{\ovl{\ID}\times \R^2 }  a\circ \phi^\tau \circ \pi (z, \xi )
\mu(dz,d\xi)=\int_{\ovl{\ID}\times \R^2 }  a \circ
\pi (z, \xi ) \mu (dz,d\xi)$ for every $a\in C^0(\W)$,
$\tau\in \R$. In other words, $\pi_* \mu$ is an invariant measure of the billiard flow.

Our main result describes finer properties of semiclassical measures $\mu$ arising from quasimodes $(u_h)$  of order $h^2$. To state it, we need to introduce some more notation.
Given $\alpha_0\in \pi \IQ\cap (-\pi/2, \pi/2)$, we will denote by $
m_{ a }^{\alpha_0, E_0}(s)$ 
the operator on $L^2_{\loc , \theta}( \R)$\footnote{The notation $L^2_{\loc, \theta}( \R )$ (resp. $L^2_{\theta}( 0,2\pi)$) is used here to emphasize that the space $L^2_{\loc}(\R)$ (resp. $L^2( 0,2\pi)$) consists in functions of the variable $\theta$.} acting by
multiplication by the function
$a \left( \Phi(s, \theta,  E_0, - E_0 \sin\alpha_0) \right)$.
If $a$ is a symmetric function (or a function on $\W$), remark
that $m_{\left\langle a\right\rangle _{\alpha_0}}^{\alpha_0, E_0}$ does
not depend on the variable $s$. For our potential $V$, the
function $\left\langle V\right\rangle _{\alpha_0}\circ\Phi$
depends only on $\theta$. 
Given $\omega\in \R/2\pi\Z$, we next define the operator
\begin{equation}
\label{e:defPetH}
P_{\alpha_0 , \omega} = -\frac{1}{2}\partial_\theta^2+
\cos^2\alpha_0\la
V\ra_{\alpha_0}\circ\Phi  , \quad \text{acting on} \quad
\cH_\omega= \{v \in L^2_{\loc}(\IR) : v(\theta + 2\pi) =
v(\theta)e^{i\omega}, \text{ for a.e. } \theta \in \R \} ,
\end{equation}
i.e. with Floquet-periodic condition. In the statements below, each
$\cH_\omega$ is identified with $L^2_{\theta}(0, 2\pi)$ by taking
restriction of functions to $(0, 2\pi)$. We are now in position to state our main result.
\begin{theorem}
\label{t:precise} Let $(u_{h})$ be a family of $O(h^2)$ quasimodes and $\mu$ be a weak-$\ast$ limit of a subsequence of $W_h$. Then, the measure $\mu$ can be decomposed into a
countable sum of nonnegative measures
$$\mu=\nu_{Leb}+\sum_{\alpha_0\in \pi\IQ\cap[-\pi/2, \pi/2]}\nu_{ \alpha_0}, \qquad \text{such that}$$

(i) Each term of the sum is carried by the
set $\{E=E_0\}$ and invariant under the billiard flow.

(ii) $\nu_{Leb}$ is of the form $\int_{|J|\leq E_0} \lambda_{E_0, J}
d\nu'(J)$ for some nonnegative measure $\nu'$ on $\R$. In
other words $\nu_{Leb}$ is a combination of Lebesgue measures on
the invariant ``tori'' $T_{(E_0, J)}$.

(iii) For $\alpha_0=\pm\frac\pi{2}$, $\nu_{\alpha_0}$ is carried by $(z, \xi)\in T^*\d\ID$, and is invariant under rotations around the origin.

(iv) For every $\alpha_0\in \pi\IQ\cap(-\pi/2, \pi/2)$, $\nu_{\alpha_0}$ is carried by the set
$\cI_{\alpha_0}\cap\{E = E_0\}$ and there exists a nonnegative measure $\ell_{\alpha_0}(d\omega)$ on $\R/2\pi\Z$, and a function
$$
\sigma_{\alpha_0}  : (\R/2\pi\Z)_\omega  \to \mathcal{L}^1_+ \big(L^2_{\theta}(0,2\pi)\big) ,
$$
integrable with respect to $\ell_{\alpha_0}$, taking values in the
set of nonnegative trace-class operators on $L^2_{\theta}(0,2\pi)$ so that
\begin{align}
\label{e:structnualpha}
\int_{\cI_{\alpha_0}} a  \, d\nu_{\alpha_0} = \int_{\cI_{\alpha_0}} \Tr_{L ^{2}_{\theta}(0, 2\pi)}\left(
m_{\la a\ra_{\alpha_0}  }^{\alpha_0, E_0}\,\sigma_{\alpha_0 }\,\right)
d \ell_{\alpha_0} , \quad \text{ for all } a\in C^\infty_c(T^*\R^2) .
\end{align}

(v) If in addition $V \in C^\infty(\overline{\ID} ; \R)$ and $(u_h)$ is a family of $o(h^2)$ quasimodes, then for $\ell_{\alpha_0}$-almost every $\omega$, we
have $\left[ P_{\alpha_0, \omega} , \sigma_{\alpha_0}(\omega) \right]
=0$ in $L^2_\theta(0,2\pi)$.
\end{theorem}

\section{Sketch of proofs}
\label{section: sketch}
\textbf{Step 1: Decomposition of an invariant measure of the billiard.}
Phase space can be partitioned into $\ovl{\ID}\times (\R^2\setminus\{0\})  =\alpha^{-1}\left(\pi\IQ\cap[-\pi/2, \pi/2]\right) \sqcup \alpha^{-1}\left(\R\setminus\pi\IQ\right)$, where $\alpha$ is the function defined in \S \ref{s:billiard}.
It follows that the invariant measure $\mu$ on $ \ovl{\ID}\times \R^2$
decomposes as a sum of nonnegative measures:%
\begin{equation}
\mu=\mu|_{ \alpha\not\in \pi\IQ}+\sum_{r\in \IQ\cap[-1/2,
1/2]}\mu|_{ \alpha=r\pi} .
\label{dec}%
\end{equation}
Since $\mu$ is a nonnegative invariant
measure on $\W$, supported in $\{|\xi|=E_0\}$, the same is true for every term in the decomposition \eqref{dec}. Moreover, $\mu|_{ \alpha\not\in
\pi\IQ}$ is invariant under the rotation flow $R^\tau$, as well
as $\mu|_{ \alpha=\pm \pi/2}$.
The assertion for $\alpha=\pm \pi/2$ comes
from the fact that the rotation flow coincides with the billiard
flow (up to time change) on the set $\{\alpha=\pm \pi/2\}$. The
assertion for $\alpha\not\in \pi\IQ$ is a standard fact: for any given value $\alpha_0$ (such that
$\alpha_0\not\in\pi\IQ$) we can find $T=T(\alpha_0)>0$ such that
$\phi^{T}$ coincides with an irrational rotation on the set
$\{\alpha=\alpha_0\}$. 
Thus, for $\alpha\not\in \pi\IQ$ or $\alpha=\pm \pi/2$, there is
nothing to prove to get Theorem \ref{t:precise}.
Hence, it only remains to study each invariant measure $\mu|_{  \alpha=\alpha_0}$, where $\alpha_0\in \pi\IQ\cap (-\pi/2, \pi/2)$ is fixed. This is the aim of the remainder of the proof.

\bigskip
\noindent
\textbf{Step 2: Second microlocalization on $\cI_{\alpha_0}$\label{s:coordIla0}.}
The angle $\alpha_0\in \pi\IQ\cap (-\pi/2, \pi/2)$ being fixed, we wish to study the concentration of $W_h$ around the set $\{J=-E\sin\alpha_0\}$. Since the limit measure $(\Phi^{-1})_* \mu$ is supported on the set $\{E =E_0\}$ this is equivalent to studying the concentration of $W_h$ around $\{J=-E_0 \sin\alpha_0\}$. 
For this, we define an appropriate class of symbols depending on
an additional variable $\eta$, which later in the calculations
will be identified with $\frac{J'}h$ for $J'= J+ E_0 \sin\alpha_0$. We denote by $\cS$ the functions $b = b(s, \theta, E, J',\eta )$ supported in $E$ away from $0$ and $+\infty$, positively homogeneous of degree zero at infinity in the variable $\eta$. We say that $b \in \cS^\sigma$ if $b \in \cS$, and $b$ and its derivatives are symmetric with respect to the boundary, which means that $ b \left(\cos\alpha,\theta, E,  J', \eta\right)= b \left(-\cos\alpha,\theta + \pi+2\alpha , E, J', \eta\right)$. 
We now introduce two auxiliary distributions which describe more precisely how $W_h$ concentrates on the set $\{E =E_0\} \cap \{J=-E_0 \sin\alpha_0\}$.
Let $\chi\in C_{c}^{\infty}\left(  \mathbb{R}\right)  $ be a
nonnegative cut-off function that is identically equal to one near
the origin and let $R>0$. For
$b\in\mathcal{S} $, we define, with $v_h = e^{-i\theta E_0 \sin\alpha_0} \scrU u_h$ and $J' = J+E_0 \sin\alpha_0$,%
$$
\left\langle w_{h,R}^{\alpha_0}   ,b\right\rangle := 
 \left\la v_h,  \Op_h\left( \left(  1-\chi\left(  \frac{ J'}{Rh}\right)\right) \chi_0(\theta)b(s, \theta,E, J', \frac{J'}{h})\right) v_h\right\ra_{L^2(\R_s\times \R_\theta)}  ,
$$
\begin{equation*}
\left\langle w_{\alpha_0, h,R}    ,b\right\rangle :=
 \left\la  v_h,  \Op_h\left(  \chi\left(  \frac{ J'}{Rh} \right) \chi_0(\theta)b(s, \theta,E, J', \frac{ J'}{h} )\right) v_h \right\ra_{L^2(\R_s\times \R_\theta )}    \label{subL} .
\end{equation*}
The Calder\'{o}n-Vaillancourt theorem ensures that both $
w_{h,R}^{\alpha_0}$ and $  w_{\alpha_0,h,R}$ are bounded in
$\mathcal{S}^{\prime}$.
After possibly extracting subsequences, we
have the existence of a limit: for every $b\in \mathcal{S} $,%
\begin{equation}
  \left\langle {\mu}^{\alpha_0
}   ,b\right\rangle :=\lim_{R\rightarrow\infty}%
\lim_{h\rightarrow0^{+}}   \left\langle w_{h,R}^{\alpha_0}
,b\right\rangle ,
\quad 
\text{ and } \quad 
   \left\langle {\mu}_{\alpha_0
}   ,b\right\rangle  :=\lim_{R\rightarrow\infty}%
\lim_{h\rightarrow0^{+}}   \left\langle
 w_{\alpha_0,h,R}   ,b\right\rangle . \label{doublelim}%
\end{equation}
These two limit distributions enjoy the following preliminary properties:

\begin{proposition}
\label{thm 1st2micro}
(i) The distribution ${\mu}^{\alpha_0}$ is a nonnegative
Radon measure. In addition, ${\mu}^{\alpha_0} $ is nonnegative,
$0$-homogeneous and supported at infinity in the variable $\eta$
; hence, ${\mu}^{\alpha_0 }$ may be
identified with a nonnegative measure on
$\R^4\times\{-1, +1\}$.

(ii) The projection of ${\mu}_{\alpha_0 } $ on $\R^4_{s, \theta, E,
J'}$, is a nonnegative measure, carried on $\{J'=0\}$, which we denote $\nu_{\alpha_0} = \int_\R {\mu}_{\alpha_0 }(d\eta)$ (in view of the statement of Theorem~\ref{t:precise}).

(iii) If $(u_h)$ is a family of $o(h)$ quasimodes, the distributions ${\mu}_{\alpha_0}  $ and ${\mu}^{\alpha_0}$ are carried by
the set $\{E=E_0 \}$ and satisfy%
\[
 \la{\mu}_{\alpha_0}  ,  \,\partial_s b \ra=0  ,\quad \la{\mu}^{\alpha_0} ,  \,\partial_s b\ra
= 0 , \quad \text{ for every }b\in\mathcal{S}^\sigma .
\]
\end{proposition}
In particular, Item~(iii) states that both ${\mu}_{\alpha_0}$ and ${\mu}^{\alpha_0}$ are, as $\mu$, supported by the set $\{|\xi|=E_0\}$ and invariant under the billiard flow.

\begin{theorem}
\label{Thm Properties}
Assume $(u_h)$ is a family of $O(h^2)$ quasimodes. 
Then, the measure $\mu^{\alpha_0}$ restricted to $\cI_{\alpha_0}$ satisfies the additional invariance property:
$\la{\mu}^{\alpha_0} |_{\cI_{\alpha_0}}  , \partial _{\theta} b \ra =0$, for every $b$ in $\cS^\sigma$.
 \end{theorem}
This is the key point to prove that, once projected to the $(s,\theta, E,J)$ variables, ${\mu}^{\alpha_0} |_{\cI_{\alpha_0}}$ is proportional to the Lebesgue measure on $\mathcal{I}_{\alpha_0}$, and hence contributes to $\nu_{Leb}$ in the statement of Theorem~\ref{t:precise}. The proof of Theorem~\ref{Thm Properties} relies on the equation \eqref{e:SR2} and involves a commutator argument. Technical problems arise when dealing with the boundary term $h^2 \d_n u_h|_{\d \ID} \otimes \delta_{\d \ID}$ : we need to go back and forth from action angles variables to polar coordinates (in which the Dirichlet boundary condition is easily expressed), developing the Fourier integral operator involved up to second order.

There remains now to study the structure of the distribution ${\mu}_{\alpha_0}$ and its invariance properties.

\bigskip
\noindent
\textbf{Step 3: Structure and propagation of ${\mu}_{\alpha_0}$.}
\begin{proposition}
There exists a nonnegative $\mathcal{L}^{1}\left( L^2_\theta(0, 2\pi) \right)$-valued measure ${\rho}_{\alpha_0}$, on $ \R/2\pi\Z_\omega \times \R_s$,  
 supported in $\{s\in [-\cos\alpha_0, \cos\alpha_0] \}$, such that
for every $b\in\cS$,
 \begin{equation}\label{e:}
\int b(s, \theta,E, J, \eta )\mu_{\alpha_0}(ds, d\theta, dE , dJ, d\eta)
=\Tr_{L ^{2}_\theta (0, 2\pi)}\int  b(s, \theta, E_0, 0, D_\theta) \,
\,\rho_{\alpha_0 } (d\omega,d s) .
\end{equation}
\end{proposition}
Similarly to Proposition \ref{thm 1st2micro} (iii), one can prove that the operator-valued measure $\rho_{\alpha_0}$ satisfies some invariance property with respect to $s$-translation. 
The very particular structure of $\mu_{\alpha_0}$ exhibited in \eqref{e:} is sufficient to prove that its projection on the variables $(s, \theta)$ is absolutely continuous.
Thus, this is also the case for the measure $\nu_{\alpha_0} = \int_\R {\mu}_{\alpha_0 }(d\eta)$ appearing in Theorem \ref{t:precise}.

The operator-valued measure $\rho_{\alpha_0}$ also possesses an additional (two-microlocal) invariance
property that we now explain. 
Setting $\overline{\rho}_{\alpha_0}(d \omega) = \int \rho_{\alpha_0}(d \omega, d s)$ and according to~\cite[Appendix]{GerardMDM91}, there exists a nonnegative measure $\ell_{\alpha_0}(d\omega)$ on $\R/2\pi\Z$, and a function $\sigma_{\alpha_0}  : (\R/2\pi\Z)_\omega  \to \mathcal{L}^1_+ \big(L^2_{\theta}(0,2\pi)\big)$, integrable with respect to $\ell_{\alpha_0}$, such that $\overline{\rho}_{\alpha_0} =\sigma_{\alpha_0} \ell_{\alpha_0}$.

\begin{theorem}\label{p:invop}
Assume that $V \in C^\infty_c(\overline \ID  ;\R)$ and that $(u_h)$ is a family of $o(h^2)$ quasimodes. Then, for $\ell_{\alpha_0}$ almost every $\omega$, we have $
\left[P_{\alpha_0, \omega} , \sigma_{\alpha_0}(\omega) \right]= 0$ in $\mathcal{H}_\omega$, where $P_{\alpha_0, \omega}$ is defined in~\eqref{e:defPetH}.
\end{theorem}
This commutation property implies that both operators are simultaneously diagonal. Combined with a unique continuation principle for eigenfunctions of the elliptic operator $P_{\alpha_0, \omega}$ from a nonempty open set, this is a key point in the proof of the observability/resolvent estimates, Corollaries \ref{t:obs-i} and~\ref{t:resolv-estim} (the paper ~\cite[Section~10]{AnL} contains a similar argument on the torus).

\footnotesize
\bibliographystyle{alpha}
\bibliography{biblio-old-and-cras}

\end{document}

%% file: def.tex
\newcommand{\nwc}{\newcommand}
\nwc{\nwt}{\newtheorem}
\nwt{coro}{Corollary}
\nwt{ex}{Example}
\nwt{prop}{Proposition}
\nwt{defin}{Definition}

\nwc{\mf}{\mathbf} 
\nwc{\blds}{\boldsymbol} 
\nwc{\ml}{\mathcal} 

\nwc{\lam}{\lambda}
\nwc{\del}{\delta}
\nwc{\Del}{\Delta}
\nwc{\Lam}{\Lambda}
\nwc{\elll}{\ell}

\nwc{\IA}{\mathbb{A}} 
\nwc{\IB}{\mathbb{B}} 
\nwc{\IC}{\mathbb{C}} 
\nwc{\ID}{\mathbb{D}} 
\nwc{\IE}{\mathbb{E}} 
\nwc{\IF}{\mathbb{F}} 
\nwc{\IG}{\mathbb{G}} 
\nwc{\IH}{\mathbb{H}} 
\nwc{\IN}{\mathbb{N}} 
\nwc{\IP}{\mathbb{P}} 
\nwc{\IQ}{\mathbb{Q}} 
\nwc{\IR}{\mathbb{R}} 
\nwc{\IS}{\mathbb{S}} 
\nwc{\IT}{\mathbb{T}} 
\nwc{\IZ}{\mathbb{Z}} 

\def\bbleft{{\mathchoice {[\mskip-3mu {[}} {[\mskip-3mu {[}}{[\mskip-4mu {[}}{[\mskip-5mu {[}}}}
\def\bbright{{\mathchoice {]\mskip-3mu {]}} {]\mskip-3mu {]}}{]\mskip-4mu {]}}{]\mskip-5mu {]}}}}
\nwc{\setK}{\bbleft 1,K \bbright}
\nwc{\setN}{\bbleft 1,\cN \bbright}
 
\nwc{\va}{{\bf a}}
\nwc{\vb}{{\bf b}}
\nwc{\vc}{{\bf c}}
\nwc{\vd}{{\bf d}}
\nwc{\ve}{{\bf e}}
\nwc{\vf}{{\bf f}}
\nwc{\vg}{{\bf g}}
\nwc{\vh}{{\bf h}}
\nwc{\vi}{{\bf i}}
\nwc{\vI}{{\bf I}}
\nwc{\vj}{{\bf j}}
\nwc{\vk}{{\bf k}}
\nwc{\vl}{{\bf l}}
\nwc{\vm}{{\bf m}}
\nwc{\vM}{{\bf M}}
\nwc{\vn}{{\bf n}}
\nwc{\vo}{{\it o}}
\nwc{\vp}{{\bf p}}
\nwc{\vq}{{\bf q}}
\nwc{\vr}{{\bf r}}
\nwc{\vs}{{\bf s}}
\nwc{\vt}{{\bf t}}
\nwc{\vu}{{\bf u}}
\nwc{\vv}{{\bf v}}
\nwc{\vw}{{\bf w}}
\nwc{\vx}{{\bf x}}
\nwc{\vy}{{\bf y}}
\nwc{\vz}{{\bf z}}
\nwc{\bal}{\blds{\alpha}}
\nwc{\bep}{\blds{\epsilon}}
\nwc{\barbep}{\overline{\blds{\epsilon}}}
\nwc{\bnu}{\blds{\nu}}
\nwc{\bmu}{\blds{\mu}}
\nwc{\bet}{\blds{\eta}}

\nwc{\bk}{\blds{k}}
\nwc{\bm}{\blds{m}}
\nwc{\bM}{\blds{M}}
\nwc{\bp}{\blds{p}}
\nwc{\bq}{\blds{q}}
\nwc{\bn}{\blds{n}}
\nwc{\bv}{\blds{v}}
\nwc{\bw}{\blds{w}}
\nwc{\bx}{\blds{x}}
\nwc{\bxi}{\blds{\xi}}
\nwc{\by}{\blds{y}}
\nwc{\bz}{\blds{z}}

\nwc{\cA}{\ml{A}}
\nwc{\cB}{\ml{B}}
\nwc{\cC}{\ml{C}}
\nwc{\cD}{\ml{D}}
\nwc{\cE}{\ml{E}}
\nwc{\cF}{\ml{F}}
\nwc{\cG}{\ml{G}}
\nwc{\cH}{\ml{H}}
\nwc{\cI}{\ml{I}}
\nwc{\cJ}{\ml{J}}
\nwc{\cK}{\ml{K}}
\nwc{\cL}{\ml{L}}
\nwc{\cM}{\ml{M}}
\nwc{\cN}{\ml{N}}
\nwc{\cO}{\ml{O}}
\nwc{\cP}{\ml{P}}
\nwc{\cQ}{\ml{Q}}
\nwc{\cR}{\ml{R}}
\nwc{\cS}{\ml{S}}
\nwc{\cT}{\ml{T}}
\nwc{\cU}{\ml{U}}
\nwc{\cV}{\ml{V}}
\nwc{\cW}{\ml{W}}
\nwc{\cX}{\ml{X}}
\nwc{\cY}{\ml{Y}}
\nwc{\cZ}{\ml{Z}}

\nwc{\fA}{\mathfrak{a}}
\nwc{\fB}{\mathfrak{b}}
\nwc{\fC}{\mathfrak{c}}
\nwc{\fD}{\mathfrak{d}}
\nwc{\fE}{\mathfrak{e}}
\nwc{\fF}{\mathfrak{f}}
\nwc{\fG}{\mathfrak{g}}
\nwc{\fH}{\mathfrak{h}}
\nwc{\fI}{\mathfrak{i}}
\nwc{\fJ}{\mathfrak{j}}
\nwc{\fK}{\mathfrak{k}}
\nwc{\fL}{\mathfrak{l}}
\nwc{\fM}{\mathfrak{m}}
\nwc{\fN}{\mathfrak{n}}
\nwc{\fO}{\mathfrak{o}}
\nwc{\fP}{\mathfrak{p}}
\nwc{\fQ}{\mathfrak{q}}
\nwc{\fR}{\mathfrak{r}}
\nwc{\fS}{\mathfrak{s}}
\nwc{\fT}{\mathfrak{t}}
\nwc{\fU}{\mathfrak{u}}
\nwc{\fV}{\mathfrak{v}}
\nwc{\fW}{\mathfrak{w}}
\nwc{\fX}{\mathfrak{x}}
\nwc{\fY}{\mathfrak{y}}
\nwc{\fZ}{\mathfrak{z}}

\nwc{\tA}{\widetilde{A}}
\nwc{\tB}{\widetilde{B}}
\nwc{\tE}{E^{\vareps}}
\nwc{\tk}{\tilde k}
\nwc{\tN}{\tilde N}
\nwc{\tP}{\widetilde{P}}
\nwc{\tQ}{\widetilde{Q}}
\nwc{\tR}{\widetilde{R}}
\nwc{\tV}{\widetilde{V}}
\nwc{\tW}{\widetilde{W}}
\nwc{\ty}{\tilde y}
\nwc{\teta}{\tilde \eta}
\nwc{\tdelta}{\tilde \delta}
\nwc{\tlambda}{\tilde \lambda}
\nwc{\ttheta}{\tilde \theta}
\nwc{\tvartheta}{\tilde \vartheta}
\nwc{\tPhi}{\widetilde \Phi}
\nwc{\tpsi}{\tilde \psi}
\nwc{\tmu}{\tilde \mu}

\nwc{\To}{\longrightarrow} 

\nwc{\ad}{\rm ad}
\nwc{\eps}{\epsilon}
\nwc{\ep}{\epsilon}
\nwc{\vareps}{\varepsilon}

\def\ep{\epsilon}

\def\Tr{{\rm Tr}}

\def\sq2{\sqrt{2}}

\def\t2{{\mathbb T}^2}
\def\s2{{\mathbb S}^2}

\def\R{\mathbb{R}}

\def\Z{\mathbb{Z}}
\def\C{\mathbb{C}}

\nwc{\lap}{\bigtriangleup}
\nwc{\rest}{\restriction}
\nwc{\Diff}{\operatorname{Diff}}
\nwc{\diam}{\operatorname{diam}}
\nwc{\Res}{\operatorname{Res}}
\nwc{\Spec}{\operatorname{Spec}}
\nwc{\Vol}{\operatorname{Vol}}
\nwc{\Op}{\operatorname{Op}}
\nwc{\supp}{\operatorname{supp}}
\nwc{\Span}{\operatorname{span}}

\nwc{\dia}{\varepsilon}
\nwc{\cut}{f}
\nwc{\qm}{u_\hbar}

\def\hto0{\xrightarrow{\hbar\to 0}}

\def\rto0{\xrightarrow{r\to 0}}

\nwc{\la}{\langle}
\nwc{\ra}{\rangle}
\nwc{\lp}{\left(}
\nwc{\rp}{\right)}

\nwc{\bequ}{\begin{equation}}
\nwc{\be}{\begin{equation}}
\nwc{\ben}{\begin{equation*}}
\nwc{\bea}{\begin{eqnarray}}
\nwc{\bean}{\begin{eqnarray*}}
\nwc{\bit}{\begin{itemize}}
\nwc{\bver}{\begin{verbatim}}

\nwc{\eequ}{\end{equation}}
\nwc{\ee}{\end{equation}}
\nwc{\een}{\end{equation*}}
\nwc{\eea}{\end{eqnarray}}
\nwc{\eean}{\end{eqnarray*}}
\nwc{\eit}{\end{itemize}}
\nwc{\ever}{\end{verbatim}}